\documentclass{article}
\usepackage{graphicx}
\usepackage{amsmath}
\usepackage{amsfonts}
\usepackage{amssymb}

\newtheorem{corollary}{Corollary}

\newtheorem{proposition}{Proposition}

\newenvironment{proof}[1][Proof]{\textbf{#1.} }{\ \rule{0.5em}{0.5em}}

\begin{document}

\title{Symplectic forms on six dimensional real solvable Lie algebras I}
\author{Rutwig Campoamor-Stursberg\\Dpto. Geometria y Topologia\\Fac. CC. Matem\'{a}ticas U.C.M.
\\Plaza de Ciencias, 3\\E-28040 Madrid, Spain\\rutwig@mat.ucm.es}
\date{}
\maketitle
\begin{abstract}
We analyze symplectic forms on six dimensional real solvable and non-nilpotent Lie algebras.
More precisely, we obtain all those algebras endowed with a symplectic form that decompose as the
direct sum of two ideals or are indecomposable solvable algebras with a four dimensional nilradical.
\end{abstract}

\section{Introduction}

Symplectic forms on Lie groups and algebras naturally appear in
the context of Poisson geometry, the study of Hamiltonians in
Mechanics and various other physical and geometrical problems,
like the cotangent bundle of differentiable manifolds. Symplectic
forms are also useful to construct other geometrical structures,
like K\"ahler manifolds, in the case of combination of symplectic
and compatible complex structures \cite{Do}. Many works have been
devoted to obtain conditions for constructing and classifying
symplectic structures on manifolds, groups or algebras, and
although no universal characterization has been obtained yet,
various general procedures of interest have been obtained
\cite{Chu,Re}. Models for Lie algebras admitting symplectic forms
have been developed in \cite{Bo}, which is related to the problem
of determining the structure of Lie groups having symplectic
forms. Various reductions have been obtained in this sense, which
simplify the problem to the analysis of solvable Lie algebras and
algebras with nontrivial Levi decomposition and a solvable
non-nilpotent radical \cite{Chu,Ha}. For fixed dimensions there
only exist complete results up to dimension four, as well as some
special cases in higher dimension \cite{AC1,C11,Co,Go1,Kr}.
Symplectic structures on four dimensional real Lie algebras have
been classified in different contexts (see e.g. \cite{Re} or
\cite{Ov1} for a recent review), and interesting applications to
the Monge-Amp\`ere equations were developed in \cite{Kr}. In
dimension six, only the nilpotent Lie algebras have been
systematically analized for symplectic forms \cite{Go1}, in
combination with additional geometrical structures (see e.g.
\cite{Co} and references therein). Results for nilpotent Lie
algebras of maximal nilindex have been obtained in \cite{Go},
while Lie algebras in dimension $n\leq 8$ endowed with exact
symplectic forms were determined in \cite{C24}. This special case
of symplectic forms is closely related to the problem of
determining the invariant functions for the coadjoint
representation of Lie groups.

In this work we begin with the systematic computation of symplectic structures on real solvable Lie
algebras of dimension six. We focus only on solvable non-nilpotent Lie algebras, basing on the classification
obtained by various authors. More precisely, we obtain the possible symplectic forms on six dimensional
solvable Lie algebras that either decompose as the direct sum of two lower dimensional ideals or are
indecomposable with a four dimensional nilradical. This covers all but one case, corresponding to
indecomposable algebras with five dimensional nilradical.

\medskip

Unless otherwise stated, any Lie algebra $\frak{g}$ considered in this work is
 defined over the field $\mathbb{R}$ of real numbers. We
convene that nonwritten brackets are either zero or obtained by antisymmetry. Abelian Lie algebras of dimension $n$ are denoted by $nL_{1}$.

\section{Symplectic structures on Lie groups}
Given a Lie algebra $\frak{g}$ with structure tensor $\left\{  C_{ij}%
^{k}\right\}  $ over a basis $\left\{  X_{1},..,X_{n}\right\}  $, the
identification of the dual space $\frak{g}^{\ast}$ with the left-invariant
Pfaffian forms on a Lie group whose algebra is isomorphic to $\frak{g}$ allows
to define an exterior differential $d$ on $\frak{g}^{\ast}$ by
\begin{equation}
d\omega\left(  X_{i},X_{j}\right)  =-C_{ij}^{k}\omega\left(  X_{k}\right)
,\;\omega\in\frak{g}^{\ast}.
\end{equation}
Therefore we can rewrite any Lie algebra $\frak{g}$ as a closed system of $2$-forms%
\begin{equation}
d\omega_{k}=-C_{ij}^{k}\omega_{i}\wedge\omega_{j},\;1\leq i<j\leq\dim\left(
\frak{g}\right)  ,
\end{equation}
called the Maurer-Cartan equations of $\frak{g}$. The closure condition $d^{2}\omega_{i}=0$ for all $i$ is equivalent to the Jacobi condition. Let $\mathcal{L}(\frak{g})=\mathbb{R}\left\{  d\omega_{i}\right\}  _{1\leq i\leq
\dim\frak{g}}$ be the linear subspace of $\bigwedge^{2}\frak{g}^{\ast}$
generated by the $2$-forms $d\omega_{i}$. It follows at once that
$\dim\mathcal{L}(\frak{g})=\dim\left(  \frak{g}\right)  $ if and only if $d\omega
_{i}\neq0$ for all $i$, that is, if $\dim\left(  \frak{g}\right)  =\dim\left[
\frak{g},\frak{g}\right]  $ holds. If $\omega=a^{i}d\omega_{i}\,\;\left(
a^{i}\in\mathbb{R}\right)  $ is an element of $\mathcal{L}(\frak{g})$, there exists a positive integer
$j_{0}\left(  \omega\right)$ such that
\begin{equation}
\bigwedge^{j_{0}\left(  \omega\right)  }\omega\neq0,\quad \bigwedge
^{j_{0}\left(  \omega\right)  +1}\omega\equiv0.
\end{equation}
This equation shows that $r\left(  \omega\right)  =2j_{0}\left(  \omega\right)
$ is the rank of the 2-form $\omega$. Define
\begin{equation}
j_{0}\left(  \frak{g}\right)  =\max\left\{  j_{0}\left(  \omega\right)
\;|\;\omega\in\mathcal{L}(\frak{g})\right\}.
\end{equation}
This quantity $j_{0}\left(  \frak{g}\right)  $ depends only on the
structure of $\frak{g}$ and is a numerical invariant of $\frak{g}$ \cite{C43}.\newline
An even dimensional Lie group $G$ is said to carry a left invariant symplectic
structure if it possesses a left invariant closed $2$-form $\omega$ of maximal
rank. At the Lie algebra level, this implies the existence of the form $\omega\in\bigwedge^{2}\frak{g}^{*}$ such  that
\begin{eqnarray}
d\omega=0\\
\bigwedge^{n}\omega\neq 0,
\end{eqnarray}
where $2n=\dim(\frak{g})$. We say that $\frak{g}$ is endowed with a symplectic structure.
For example, any Lie algebra in dimension 2 has a symplectic structure.
For the abelian algebra the assertion is trivial, while for the
affine Lie algebra $\frak{r}_{2}=\left\{  X_{1},X_{2}|\quad[X_{1},X_{2}%
]=X_{2}\right\}  $ we have $\omega\in\frak{r}_{2}^{*}\wedge\frak{r}_{2}^{*}$ defined
by
\begin{equation}
\omega=\omega_{1}\wedge\omega_{2}.
\end{equation}
This form is closed and of maximal rank. In particular, $\omega\in\mathcal{L}(\frak{r}_{2})$.
Symplectic form of this special kind are called exact symplectic forms, and they are of interest in
the analysis of invariants of Lie algebras \cite{AC1,C43}. That is, a Lie algebra is exact symplectic if
it is symplectic and the form is
moreover exact.

The structure of a Lie algebra plays an essential role in the
existence of such forms. We briefly recall some results which will be used
in this work.

\begin{proposition}
{\rm \cite{Chu}} Let $\frak{g}$ be a Lie algebra. Then following conditions hold:

\begin{enumerate}
\item  If $Trace[ad(X)]=0\quad\forall
X\in\frak{g}$) and $\frak{g}$ is symplectic, then $\frak{g}$ is solvable.

\item No semisimple Lie algebra carries a symplectic form.

\item  Direct sums of semisimple and solvable Lie algebras cannot be symplectic.

\end{enumerate}
\end{proposition}

This implies that an indecomposable symplectic Lie algebra is either solvable or the semidirect product of a semisimple Lie algebra and a solvable non-nilpotent Lie algebra. In particular, any symplectic Lie algebra in dimension four is solvable. We remark that in this dimension any nilpotent Lie algebra is endowed with a symplectic form \cite{Ha}.

\subsection{Six dimensional Lie algebras}

Real Lie algebras of dimension six have been fully classified by Morozov (nilpotent real Lie algebras), Mubarakzyanov (decomposition conditions and solvable algebras with five dimensional nilradical) and Turkowski (solvable algebras with four dimensional nilradical and algebras with nontrivial Levi subalgebra) (see e.g. \cite{Mu1,Mu2,Tu} and references therein). According to the list in \cite{Tu2}, there are four indecomposable Lie algebras with nonzero Levi subalgebra, corresponding to the semidirect product of the abelian Lie algebra $3L_{1}$ with $\frak{sl}(2,\mathbb{R})$ and $\frak{so}(3)$, the semidirect product of $\frak{sl}(2,\mathbb{R})$ and the Heisenberg algebra in dimension 3 and the Lie algebra $\frak{sl}(2,\mathbb{R})\overrightarrow{\oplus}_{D_{\frac{1}{2}}\oplus D_{0}}A_{3,3}$, which is the only to possess a (exact) symplectic form. The decomposable case follows from \cite{Chu}. It remains to analyze the solvable case. Nilpotent algebras endowed with symplectic forms can be found in \cite{Co}, for which reason we do not reproduce them here. As known, the maximal nilpotent ideal (nilradical) $NR$ of a solvable Lie algebra $\frak{r}$ satisfies the inequality
\begin{equation}
\dim (NR)\leq \frac{\frak{r}}{2}. \label{FM}
\end{equation}
Therefore a solvable algebra in dimension six has a nilradical of dimensions four or five if it is indecomposable, or it is the direct sum of ideals.

\section{Decomposable solvable Lie algebras}

In this section we analize the symplectic solvable Lie algebras $\frak{r}$ that decompose as the direct sum $\frak{r}=\frak{r}_{1}\oplus\frak{r}_{2}$ of lower dimensional ideals. Since any solvable Lie algebra is supposed to be non-nilpotent, at least one of the ideals $\frak{r}_{i}$ must be non-nilpotent.

Let $\frak{r}_{1}$ and $\frak{r}_{2}$ be solvable Lie algebras of odd
dimension and let $\left\{  \omega_{1},..,\omega_{2n+1}\right\}  $ and
$\left\{  \omega_{1}^{\prime},..,\omega_{2m+1}^{\prime}\right\}  $ be bases of
$\frak{r}_{1}^{\ast}$, respectively $\frak{r}_{2}^{\ast}$. Solvability implies
that $\frak{r}_{i}\neq\left[  \frak{r}_{i},\frak{r}_{i}\right]  $ for $i=1,2$,
so that without loss of generality we can suppose that
\begin{equation}
d\omega_{1}=d\omega_{1}^{\prime}=0.
\end{equation}
In these conditions, we give a sufficiency criterion for the existence of
symplectic forms on the direct sum algebra $\frak{r}_{1}\oplus\frak{r}_{2}$:

\begin{proposition}
Suppose that there exist 2-forms $\theta=\sum_{i<j}\alpha^{ij}\omega_{i}%
\wedge\omega_{j}\in\bigwedge^{2}\frak{r}_{1}^{\ast}$ and $\theta^{\prime}%
=\sum_{k<l}\beta^{kl}\omega_{k}^{\prime}\wedge\omega_{l}^{\prime}\in
\bigwedge^{2}\frak{r}_{2}^{\ast}$ satisfying

\begin{enumerate}
\item $\alpha^{1j}=\beta^{1l}=0,\;j,l\geq2.$

\item $d\theta=d\theta^{\prime}=0,$

\item $\bigwedge^{n}\theta\neq0$ and $\bigwedge^{m}\theta^{\prime}\neq0.$
\end{enumerate}
Then $\frak{r}_{1}\oplus\frak{r}_{2}$ is endowed with a sympletic form
$\eta=\theta+\theta^{\prime}+\omega_{1}\wedge\omega_{1}^{\prime}$
\end{proposition}

\begin{proof}
Since $\alpha^{1j}=0$ for any $j$, the 2-form $\theta$ can be written as
\begin{equation}
\theta=\sum_{2\leq i<j\leq2n+1}\alpha^{ij}\omega_{i}\wedge\omega_{j}.
\end{equation}
By assumption, the $n^{th}$-exterior product $\bigwedge^{n}\theta$ is not
zero, so that reordering the basis if necessary, we can suppose that
\begin{equation}
\prod_{k=1}^{n}\alpha^{2k,2k+1}\neq0.
\end{equation}
In consequence we obtain
\begin{equation}
\bigwedge^{n}\theta=\left(  n!\prod_{k=1}^{n}\alpha^{2k,2k+1}\right)
\omega_{2}\wedge\omega_{3}\wedge...\wedge\omega_{2n}\wedge\omega_{2n+1}.
\end{equation}
A similar expression holds for $\theta^{\prime}$. Clearly the 2-form
\begin{equation}
\eta=\theta+\theta^{\prime}+\omega_{1}\wedge\omega_{1}^{\prime}%
\end{equation}
belongs to $\left(  \frak{r}_{1}\oplus\frak{r}_{2}\right)  ^{\ast}$. Further%
\begin{equation}
d\eta=d\theta+d\theta^{\prime}+d\omega_{1}\wedge\omega_{1}^{\prime}-\omega
_{1}\wedge d\omega_{1}^{\prime}=0
\end{equation}
by condition 2, showing that $\eta$ is closed. Finally,
\begin{equation}
\bigwedge^{n+m+1}\eta=\left(  \left(  n+m+1\right)  !\prod_{k=1}^{n}%
\alpha^{2k,2k+1}\prod_{l=1}^{m}\right)  \omega_{1}\wedge...\wedge\omega
_{2n+1}\wedge\omega_{1}^{\prime}\wedge...\wedge\omega_{2m+1}\neq0,
\end{equation}
showing that $\eta$ is of maximal rank.
\end{proof}

The result has an interesting consequence concerning odd-dimensional solvable Lie algebras.
A 1-form $\omega\in\frak{r}^{*}$ is called a linear contact form if
\begin{equation}
\omega \wedge \left( d\omega _{\mu }\right) ^{n}\neq 0.
\end{equation}
In particular, the left invariant Pfaff form induced by $\omega$ over the Lie group having $\frak{r}$ as Lie algebra is a contact form in the classical sense.

\begin{corollary}
Let $\frak{r}$ be solvable Lie algebra and $\left\{  \omega_{1},..,\omega
_{2n+1}\right\}  $ be a basis of $\frak{r}^{\ast}$. Suppose that $d\omega
_{1}=0.$ and that there exists  a 2-form $\theta=\sum_{i<j}\alpha^{ij}%
\omega_{i}\wedge\omega_{j}\in\bigwedge^{2}\frak{r}^{\ast}$ such that

\begin{enumerate}
\item $\alpha^{1j}=0,\;j\geq2.$

\item $d\theta=0,$

\item $\bigwedge^{n}\theta\neq0$
\end{enumerate}

If $\theta\in\mathcal{L}\left(  \frak{r}\right)  $, then $\frak{r}$ is endowed
with a linear contact form.
\end{corollary}

\begin{proof}
By assumption the 2-form $\theta$ can be written as
\begin{equation}
\theta=\sum_{2\leq i<j\leq2n+1}\alpha^{ij}\omega_{i}\wedge\omega_{j}.
\end{equation}
Without loss of generality we can suppose that
\begin{equation}
\bigwedge^{n}\theta=\left(  n!\prod_{k=1}^{n}\alpha^{2k,2k+1}\right)
\omega_{2}\wedge\omega_{3}\wedge...\wedge\omega_{2n}\wedge\omega_{2n+1}.
\end{equation}
If $\theta\in\mathcal{L}\left(  \frak{r}\right)  $, then there exist scalars
$a_{i_{1}},..,a_{i_{k}}$ such that
\begin{equation}
\theta=a_{i_{1}}d\omega_{i_{1}}+..+a_{i_{k}}d\omega_{i_{k}}.
\end{equation}
Further $i_{j}\neq1$ for $j\in\left\{  1,..,k\right\}  $ since $d\omega_{1}%
=0$. Define the linear form
\begin{equation}
\eta=\omega_{1}+a_{i_{1}}\omega_{i_{1}}+..+a_{i_{k}}\omega_{i_{k}}%
\end{equation}
Then
\begin{eqnarray}
\eta\wedge\left(  \bigwedge^{n}d\eta\right)  =\left(  \omega_{1}+a_{i_{1}%
}\omega_{i_{1}}+..+a_{i_{k}}\omega_{i_{k}}\right)  \wedge\left(  \bigwedge
^{n}\theta\right)  =\nonumber\\
\left(  n!\prod_{k=1}^{n}\alpha^{2k,2k+1}\right)
\omega_{1}\wedge\omega_{2}\wedge\omega_{3}\wedge...\wedge\omega_{2n}%
\wedge\omega_{2n+1}\neq0,
\end{eqnarray}
thus $\eta$ is a contact form.
\end{proof}

\subsection{$\dim \frak{r}_{1}=\dim \frak{r}_{2}=3$}

There are five isomorphism classes (two of them depending on parameters) of indecomposable solvable Lie algebras in dimension 3, whose Maurer-Cartan equations are given in Table 4 of the appendix. From these algebras, only one is nilpotent (the Heisenberg algebra). We have therefore 14 solvable non-nilpotent algebras $\frak{r}$ which decompose as the direct sum of two three dimensional ideals. To determine the possible symplectic forms on these algebras, we use proposition 2. The resulting algebras admitting such structures are given in table 1.

\begin{table}[h]
\caption{Decomposable algebras }
\begin{tabular}
[c]{llll}%
Algebra & exact & Symplectic form & Condition \\\hline
$A_{3,1}\oplus A_{3,4}^{-1}$ & no & $a_{1}\omega_{1}\wedge\omega_{2}%
+a_{2}\omega_{1}\wedge\omega_{3}+a_{3}\omega_{2}\wedge\omega_{3}+a_{4}%
\omega_{3}\wedge\omega_{5}+$ & \\
&  & $a_{5}\omega_{3}\wedge\omega_{6}+a_{6}\omega_{4}\wedge\omega_{5}%
+a_{7}\omega_{4}\wedge\omega_{6}+a_{8}\omega_{5}\wedge\omega_{6}$ &
$a_{1}\left(  a_{4}a_{7}-a_{5}a_{6}\right)  \neq0$\\
$A_{3,1}\oplus A_{3,5}^{0}$ & no & $a_{1}\omega_{1}\wedge\omega_{2}%
+a_{2}\omega_{1}\wedge\omega_{3}+a_{3}\omega_{2}\wedge\omega_{3}+a_{4}%
\omega_{3}\wedge\omega_{5}+$ & \\
&  & $a_{5}\omega_{3}\wedge\omega_{6}+a_{6}\omega_{4}\wedge\omega_{5}%
+a_{7}\omega_{4}\wedge\omega_{6}+a_{8}\omega_{5}\wedge\omega_{6}$ &
$a_{1}\left(  a_{4}a_{7}-a_{5}a_{6}\right)  \neq0$\\
$A_{3,4}^{-1}\oplus A_{3,4}^{-1}$ & no & $a_{1}\omega_{1}\wedge\omega
_{2}+a_{2}\omega_{1}\wedge\omega_{3}+a_{3}\omega_{2}\wedge\omega_{3}%
+a_{4}\omega_{3}\wedge\omega_{6}+$ & \\
&  & $+a_{5}\omega_{4}\wedge\omega_{5}+a_{6}\omega_{4}\wedge\omega_{6}%
+a_{7}\omega_{5}\wedge\omega_{6}$ & $a_{1}a_{4}a_{5}\neq0$\\
$A_{3,4}^{-1}\oplus A_{3,5}^{0}$ & no & $a_{1}\omega_{1}\wedge\omega_{2}%
+a_{2}\omega_{1}\wedge\omega_{3}+a_{3}\omega_{2}\wedge\omega_{3}+a_{4}%
\omega_{3}\wedge\omega_{6}+$ & \\
&  & $+a_{5}\omega_{4}\wedge\omega_{5}+a_{6}\omega_{4}\wedge\omega_{6}%
+a_{7}\omega_{5}\wedge\omega_{6}$ & $a_{1}a_{4}a_{5}\neq0$\\
$A_{3,5}^{0}\oplus A_{3,5}^{0}$ & no & $a_{1}\omega_{1}\wedge\omega_{2}%
+a_{2}\omega_{1}\wedge\omega_{3}+a_{3}\omega_{2}\wedge\omega_{3}+a_{4}%
\omega_{3}\wedge\omega_{6}+$ & \\
&  & $+a_{5}\omega_{4}\wedge\omega_{5}+a_{6}\omega_{4}\wedge\omega_{6}%
+a_{7}\omega_{5}\wedge\omega_{6}$ & $a_{1}a_{4}a_{5}\neq0$\\\hline
\end{tabular}
\end{table}

\subsection{$\dim \frak{r}_{4}=\dim \frak{r}_{2}=2$}

The case of direct sums $\frak{r}=\frak{r}_{1}\oplus\frak{r}_{2}$ with $\dim \frak{r}_{4}=\dim \frak{r}_{2}=2$ follows at once from the classification of symplectic structures on four dimensional real Lie algebras carried out in \cite{Ov1}. Indeed, since any two dimensional Lie algebra has a symplectic form, the sum with any four dimensional algebra having also a symplectic form gives a six dimensional algebra. On the other hand, it is immediate that no direct sum of a four dimensional Lie algebra with no symplectic form and a two dimensional algebra can result in a six dimensional symplectic Lie algebra. Therefore the result is obtained combining the results of \cite{Ov1} with those of dimension two.

\subsection{$\dim \frak{r}_{5}=\dim \frak{r}_{2}=1$}

The Maurer-Cartan equations of the indecomposable real solvable Lie algebras in dimension five are given in tables 5 and 6 of the appendix. There are 33 cases to be analyzed. By (\ref{FM}), the nilradical of such algebras $\frak{r}$ have dimension three or four. It is not difficult to see that a direct sum $\frak{r}\oplus L_{1}$ endowed with a symplectic form must satisfy the requirements of proposition 2 (otherwise the closure of the 2-form would be violated). In particular we have $\theta^{\prime}=0$, since $\dim L_{1}=1$. The algebras admitting a symplectic form are listed in table 2.

\begin{table}[h]
\caption{$\dim \frak{r}_{5}=\dim \frak{r}_{2}=1$ }
\begin{tabular}
[c]{lll}%
Algebra & Symplectic form & Conditions\\\hline
$\frak{g}_{5,7}^{\alpha,-\alpha,-1}\oplus L_{1}$ & $\sum_{i=1}^{4}a_{i}%
d\omega_{i}+a_{5}\omega_{1}\wedge\omega_{4}+a_{6}\omega_{2}\wedge\omega
_{3}+a_{7}\omega_{5}\wedge\omega_{6}$ & $a_{5}a_{6}a_{7}\neq0$\\
$\frak{g}_{5,7}^{1,-1,-1}\oplus L_{1}$ & $\sum_{i=1}^{4}a_{i}d\omega_{i}%
+a_{5}\omega_{1}\wedge\omega_{4}+a_{6}\omega_{2}\wedge\omega_{3}+a_{7}%
\omega_{5}\wedge\omega_{6}+$ & \\
& $+a_{8}\omega_{1}\wedge\omega_{3}+a_{9}\omega_{2}\wedge\omega_{4}$ &
$a_{7}\left(  a_{5}a_{6}-a_{8}a_{9}\right)  \neq0$\\
$\frak{g}_{5,8}^{-1}\oplus L_{1}$ & $\sum_{i=1}^{4}a_{i}d\omega_{i}%
+a_{5}\omega_{1}\wedge\omega_{2}+a_{6}\omega_{1}\wedge\omega_{5}+a_{7}%
\omega_{2}\wedge\omega_{6}+$ & \\
& $+a_{8}\omega_{3}\wedge\omega_{4}+a_{9}\omega_{5}\wedge\omega_{6}$ &
$a_{8}\left(  a_{5}a_{9}-a_{6}a_{7}\right)  \neq0$\\
$\frak{g}_{5,13}^{-1,0,s}\oplus L_{1}$ & $\sum_{i=1}^{4}a_{i}d\omega_{i}%
+a_{5}\omega_{1}\wedge\omega_{2}+a_{6}\omega_{3}\wedge\omega_{4}+a_{7}%
\omega_{5}\wedge\omega_{6}$ & $a_{5}a_{6}a_{7}\neq0$\\
$\frak{g}_{5,15}^{-1}\oplus L_{1}$ & $\sum_{i=1}^{4}a_{i}d\omega_{i}%
+a_{5}\left(  \omega_{1}\wedge\omega_{4}-\omega_{2}\wedge\omega_{3}\right)
+a_{7}\omega_{2}\wedge\omega_{4}+$ & \\
& $+a_{8}\omega_{5}\wedge\omega_{6}$ & $a_{5}a_{8}\neq0$\\
$\frak{g}_{5,17}^{p,-p,\pm1}$ & $\sum_{i=1}^{4}a_{i}\omega_{i}\wedge\omega
_{5}+a_{5}\omega_{5}\wedge\omega_{6}+$ & \\
& $+a_{6}\left(  s\omega_{1}\wedge\omega_{3}+\omega_{2}\wedge\omega
_{4}\right)  +a_{7}\left(  s\omega_{1}\wedge\omega_{4}-\omega_{2}\wedge
\omega_{3}\right)  $ & $a_{5}\left(  a_{6}^{2}+a_{7}^{2}\right)  \neq0$\\
$\frak{g}_{5,17}^{0,0,s}$ & $\sum_{i=1}^{4}a_{i}d\omega_{i}+a_{5}\omega
_{1}\wedge\omega_{2}+a_{6}\omega_{3}\wedge\omega_{4}+a_{7}\omega_{5}%
\wedge\omega_{6}$ & \\
& $+a_{8}\left(  s\omega_{1}\wedge\omega_{3}+\omega_{2}\wedge\omega
_{4}\right)  +a_{9}\left(  \omega_{1}\wedge\omega_{4}-s\omega_{2}\wedge
\omega_{3}\right)  $ & $a_{7}\left(  a_{5}a_{6}-s\left(  a_{8}^{2}+a_{9}%
^{2}\right)  \right)  \neq0$\\
$\frak{g}_{5,18}^{0}\oplus L_{1}$ & $\sum_{i=1}^{4}a_{i}\omega_{i}\wedge
\omega_{5}+a_{5}\omega_{5}\wedge\omega_{6}+a_{6}\left(  \omega_{1}\wedge
\omega_{3}+\omega_{2}\wedge\omega_{4}\right)  +$ & \\
& $+a_{7}\omega_{3}\wedge\omega_{4}$ & $a_{5}a_{6}\neq0$\\
$\frak{g}_{5,19}^{\frac{-1}{2},-1}\oplus L_{1}$ & $\sum_{i=1}^{4}a_{i}%
d\omega_{i}+a_{5}\omega_{1}\wedge\omega_{3}+a_{6}\omega_{2}\wedge\omega
_{4}+a_{7}\omega_{5}\wedge\omega_{6}$ & $a_{5}a_{6}a_{7}\neq0$\\
$\frak{g}_{5,19}^{-2,2}\oplus L_{1}$ & $\sum_{i=1}^{4}a_{i}d\omega_{i}%
+a_{5}\omega_{1}\wedge\omega_{2}+a_{6}\omega_{3}\wedge\omega_{4}+a_{7}%
\omega_{5}\wedge\omega_{6}$ & $a_{5}a_{6}a_{7}\neq0$\\
$\frak{g}_{5,30}^{0}\oplus L_{1}$ & $\sum_{i=1}^{4}a_{i}d\omega_{i}%
+a_{5}\omega_{3}\wedge\omega_{5}+a_{6}\omega_{3}\wedge\omega_{6}+a_{7}%
\omega_{5}\wedge\omega_{6}$ & $a_{1}a_{6}\neq0$\\
$\frak{g}_{5,33}^{-1,0}\oplus L_{1}$ & $\sum_{i=1}^{3}a_{i}d\omega_{i}%
+a_{4}\omega_{1}\wedge\omega_{3}+a_{5}\omega_{4}\wedge\omega_{5}+a_{6}%
\omega_{4}\wedge\omega_{6}+a_{7}\omega_{5}\wedge\omega_{6}$ & $a_{2}a_{4}%
a_{6}\neq0$\\
$\frak{g}_{5,33}^{0,-1}\oplus L_{1}$ & $\sum_{i=1}^{3}a_{i}d\omega_{i}%
+a_{4}\omega_{2}\wedge\omega_{3}+a_{5}\omega_{4}\wedge\omega_{5}+a_{6}%
\omega_{4}\wedge\omega_{6}+a_{7}\omega_{5}\wedge\omega_{6}$ & $a_{1}a_{4}%
a_{7}\neq0$\\
$\frak{g}_{5,36}\oplus L_{1}$ & $\sum_{i=1}^{3}a_{i}d\omega_{i}+a_{4}%
\omega_{4}\wedge\omega_{5}+a_{5}\omega_{4}\wedge\omega_{6}+a_{6}\omega
_{5}\wedge\omega_{6}$ & $a_{1}a_{6}\neq0$\\
$\frak{g}_{5,37}\oplus L_{1}$ & $\sum_{i=1}^{3}a_{i}d\omega_{i}+a_{4}%
\omega_{4}\wedge\omega_{5}+a_{5}\omega_{4}\wedge\omega_{6}+a_{6}\omega
_{5}\wedge\omega_{6}$ & $a_{1}a_{6}\neq0$\\
&  & \\\hline
\end{tabular}
\end{table}

\section{Solvable Lie algebras with four dimensional nilradical}

As follows from the classification in \cite{Tu}, there are 33 indecomposable solvable non-nilpotent real Lie algebras in dimension six with a four dimensional nilradical. The corresponding Maurer-Cartan equations for these algebras are listed in tables 7 and 8 of the appendix. The search for symplectic structures in this case must be developed case by case, since the algebras are indecomposable.

\begin{proposition}
Let $\frak{r}$ be an indecomposable solvable non-nilpotent real Lie algebra of
dimension six. Then $\frak{r}$ is endowed with a symplectic form $\omega$ if
and only if it is isomorphic to one of the Lie algebras in table 3.
\end{proposition}

\begin{proof}
We give the detailed proof for $N_{6,1}^{\alpha\beta\gamma\delta}$, all the remaining cases are treated in a similar way.

The brackets over the basis $\left\{  N_{1},..,N_{4},X_{1},X_{2}\right\}  $
are given by
\begin{eqnarray}
\left[  X_{1},N_{1}\right]   &  =\alpha N_{1},\;\left[  X_{1},N_{2}\right]
=\gamma N_{2},\;\left[  X_{1},N_{4}\right]  =N_{4}\\
\left[  X_{2},N_{1}\right]   &  =\beta N_{1},\;\left[  X_{2},N_{2}\right]
=\delta N_{2},\;\left[  X_{2},N_{3}\right]  =N_{3},
\end{eqnarray}
where $\alpha,\beta,\gamma,\delta\in\mathbb{R}$ satisfy the restrictions $\alpha\beta\neq 0$ and $\gamma^{2}+\delta^{2}\neq 0$. Taking the dual basis $\left\{  \eta_{1},..,\eta_{4},\omega_{1},\omega
_{2}\right\}  $, the Maurer-Cartan equations of the algebra are%
\[%
\begin{tabular}
[c]{ll}%
$d\eta_{1}=\alpha\omega_{1}\wedge\eta_{1}+\beta\omega_{2}\wedge\eta_{1},$ &
$d\eta_{2}=\gamma\omega_{1}\wedge\eta_{2}+\delta\omega_{2}\wedge\eta_{2},$\\
$d\eta_{3}=\omega_{2}\wedge\eta_{3},$ & $d\eta_{4}=\omega_{1}\wedge\eta_{4}%
,$\\
$d\omega_{1}=0,$ & $d\omega_{2}=0.$%
\end{tabular}
\]
Now define an element $\omega\in\bigwedge^{2}\left(  N_{6,1}%
^{\alpha\beta\gamma\delta}\right)  ^{\ast}$ in general position by
\begin{equation}
\omega=a^{ij}\eta_{i}\wedge\eta_{j}+b^{ik}\eta_{i}\wedge\omega_{k}%
+c^{12}\omega_{1}\wedge\omega_{2},\;
\end{equation}
where $1\leq i<j\leq4,\;k=1,2$ and $a^{ij},b^{ik},c^{12}\in\mathbb{R}$. If we
impose the closure%
\begin{equation}
d\omega=a^{ij}d\eta_{i}\wedge\eta_{j}+b^{ik}d\eta_{i}\wedge\omega_{k}%
-a^{ij}\eta_{i}\wedge d\eta_{j}=0
\end{equation}
and take into account the relations satisfied by the parameters, then
the following coefficients vanish independently of their value
$\alpha,\beta,\gamma$ and $\delta:$
\begin{equation}
a^{13}=a^{14}=a^{34}=b^{31}=b^{42}=0,
\end{equation}
and since $\alpha\beta\neq0$,
\begin{equation}
b^{12}=\frac{\beta}{\alpha}b^{11}.
\end{equation}
We further obtain the following expression for the differential:%
\begin{eqnarray}
d\omega=a^{24}\delta\omega_{2}\wedge\eta_{2}\wedge\eta_{4}+a^{23}\left(
1+\delta\right)  \omega_{2}\wedge\eta_{2}\wedge\eta_{3}+a^{12}\left(
\alpha+\gamma\right)  \omega_{1}\wedge\eta_{1}\wedge\eta_{2}+\nonumber \\
\left(  \gamma b^{22}-\delta b^{21}\right)  \omega_{1}\wedge\eta_{2}%
\wedge\omega_{2}+a^{24}\left(  1+\gamma\right)  \omega_{1}\wedge\eta_{2}%
\wedge\eta_{4}+a^{23}\gamma\omega_{1}\wedge\eta_{2}\wedge\eta_{3}+\nonumber \\
+a^{12}\left(  \beta+\delta\right)  \omega_{2}\wedge\eta_{1}\wedge\eta_{2}=0. \label{GL}
\end{eqnarray}

At this stage, the analysis must be divided into several steps, according to the different
possibilities depending on the four parameters:

\begin{enumerate}
\item  Let $\gamma\neq0$. Then the closure (\ref{GL}) implies
\begin{equation}
\begin{tabular}
[c]{ll}%
$a^{23}=0,$ & $b^{22}=\frac{\delta}{\gamma}b^{21},$\\
$\delta a^{24}=0,$ & $a^{12}\left(  \alpha+\gamma\right)  =0,$\\
$a^{24}\left(  1+\gamma\right)  =0,$ & $a^{12}\left(  \beta+\delta\right)
=0.$%
\end{tabular}
\end{equation}
and we obtain the wedge product
\begin{equation}
\bigwedge^{3}\omega=6b^{32}\left(  a^{12}b^{41}+b^{11}a^{24}\right)  \eta
_{1}\wedge..\wedge\eta_{4}\wedge\omega_{1}\wedge\omega_{2}%
\end{equation}

\begin{enumerate}
\item  If $\delta=0$, then $a^{12}=0$ since $\beta$ is nonzero, and the wedge product $\bigwedge
^{3}\omega$ is nonzero if and only if
\begin{equation}
1+\gamma=0.
\end{equation}
Therefore only the Lie algebra $N_{6,1}^{\alpha,\beta,-1,0}$ has a symplectic
form, given by
\begin{eqnarray}
\omega = b^{11}\eta_{1}\wedge\omega_{1}+\frac{\beta}{\alpha}b^{11}\eta
_{1}\wedge\omega_{2}+a^{24}\eta_{2}\wedge\eta_{4}+b^{21}\eta_{2}\wedge
\omega_{1}+b^{32}\eta_{3}\wedge\omega_{2}+b^{41}\eta_{4}\wedge\omega
_{1}+\nonumber \\
+c^{12}\omega_{1}\wedge\omega_{2}= -\frac{b^{11}}{\alpha}d\eta_{1}-b^{21}d\eta_{2}-b^{32}d\eta_{3}%
-b^{41}d\eta_{4}+a^{24}\eta_{2}\wedge\eta_{4}+c^{12}\omega_{1}\wedge\omega
_{2}.
\end{eqnarray}

\item  If $\delta\neq0$, then by (\ref{GL}) we have $a^{24}=0$. We obtain that
\begin{equation}
\bigwedge^{3}\omega=6b^{32}\left(  a^{12}b^{41}\right)  \eta_{1}\wedge
..\wedge\eta_{4}\wedge\omega_{1}\wedge\omega_{2}\neq0.
\end{equation}
This product is different from zero if and only if
\begin{equation}
\alpha+\gamma=\beta+\delta=0.
\end{equation}
In this case the Lie algebra $N_{6,1}^{\alpha,\beta,-\alpha,-\beta}$ has the
symplectic form%
\begin{align*}
\omega  =a^{13}\eta_{1}\wedge\eta_{2}+b^{11}\left(  \eta_{1}\wedge\omega
_{1}+\frac{\beta}{\alpha}\eta_{1}\wedge\omega_{2}\right)  +b^{21}\left(
\eta_{2}\wedge\omega_{1}+\frac{\beta}{\alpha}\eta_{2}\wedge\omega_{2}\right)
+b^{32}\eta_{3}\wedge\omega_{2}+\\
+b^{41}\eta_{4}\wedge\omega_{1}+c^{12}%
\omega_{1}\wedge\omega_{2}  =a^{12}\eta_{1}\wedge\eta_{2}-\frac{b^{11}}{\alpha}d\eta_{1}-\frac{b^{21}%
}{\alpha}d\eta_{2}-b^{32}d\eta_{3}-b^{41}d\eta_{4}+c^{12}\omega_{1}\wedge\omega_{2}.
\end{align*}
\end{enumerate}

\item  Let $\gamma=0:$ then $\delta\neq0$ and the closure implies
\begin{align*}
a^{12} &  =0\\
a^{24} &  =0,\\
a^{23}\left(  1+\delta\right)   &  =0.
\end{align*}
In addition
\begin{equation}
\bigwedge^{3}\omega=6\frac{\beta}{\alpha}b^{11}a^{23}b^{41}\eta_{1}%
\wedge..\wedge\eta_{4}\wedge\omega_{1}\wedge\omega_{2}.
\end{equation}
It is nonzero if and only if $1+\delta=0$, and in this case the Lie algebra
$N_{6,1}^{\alpha,\beta,0,-1}$ has the symplectic form
\begin{align*}
\omega &  =b^{11}\eta_{1}\wedge\omega_{1}+\frac{\beta}{\alpha}b^{11}\eta
_{1}\wedge\omega_{2}+a^{23}\eta_{2}\wedge\eta_{3}+b^{22}\eta_{2}\wedge
\omega_{2}+b^{32}\eta_{3}\wedge\omega_{2}+b^{41}\eta_{4}\wedge\omega
_{1}+c^{12}\omega_{1}\wedge\omega_{2}\\
&  =-\frac{b^{11}}{\alpha}d\eta_{1}-b^{22}d\eta_{2}-b^{32}d\eta_{3}%
-b^{41}d\eta_{4}+a^{23}\eta_{2}\wedge\eta_{3}+c^{12}\omega_{1}\wedge\omega
_{2}.
\end{align*}
\end{enumerate}
Resuming, the Lie algebras  $N_{6,1}^{\alpha,\beta,\gamma,\delta}$ admit symplectic (and non-exact) forms if
\begin{equation*}
(\gamma,\delta)\in\left\{(0,-,1),(-1,0),(-\alpha,-\beta)\right\}.
\end{equation*}
The parametrs $\alpha$ and $\beta$ are not subjected to further constraints.
\end{proof}

\begin{table}
\caption{Indecomposable algebras with four dimensional nilradical.}
\begin{tabular}
[c]{lllll}%
Algebra & Parameters & exact & Symplectic form & Conditions\\\hline
$N_{6,1}^{\alpha,\beta,-\alpha,-\beta}$ &  & no & $\sum_{i=1}^{4}a_{i}%
d\eta_{i}+a_{5}\eta_{2}\wedge\eta_{3}+a_{6}\omega_{1}\wedge\omega_{2}$ &
$a_{3}a_{4}a_{5}\neq0$\\
$N_{6,1}^{\alpha,\beta,-1,0}$ &  & no & $\sum_{i=1}^{4}a_{i}d\eta_{i}%
+a_{5}\eta_{1}\wedge\eta_{2}+a_{6}\omega_{1}\wedge\omega_{2}$ & $a_{3}%
a_{4}a_{5}\neq0$\\
$N_{6,1}^{\alpha,\beta,0,-1}$ &  & no & $\sum_{i=1}^{4}a_{i}d\eta_{i}%
+a_{5}\eta_{2}\wedge\eta_{3}+a_{6}\omega_{1}\wedge\omega_{2}$ & $a_{1}%
a_{4}a_{5}\neq0$\\
$N_{6,2}^{0,-1,\gamma}$ &  & no & $\sum_{i=1}^{4}a_{i}d\eta_{i}+a_{5}\eta
_{2}\wedge\eta_{3}+a_{6}\omega_{1}\wedge\omega_{2}$ & $a_{1}a_{2}a_{5}\neq0$\\
$N_{6,7}^{0,\beta,0,\delta}$ &  & no & $\sum_{i=1}^{4}a_{i}d\eta_{i}+a_{5}%
\eta_{1}\wedge\eta_{2}+a_{6}\omega_{1}\wedge\omega_{2}$ & $a_{4}a_{5}\neq0$\\
$N_{6,13}^{\alpha\beta\gamma\delta}$ & $\alpha+\gamma=0,\;$ & no & $\sum
_{i=1}^{4}a_{i}d\eta_{i}+a_{5}\eta_{1}\wedge\eta_{2}+a_{6}\omega_{1}%
\wedge\omega_{2}$ & \\
& $\beta+\delta=0$ &  &  & $a_{5}\left(  a_{3}^{2}+a_{4}^{2}\right)  \neq0$\\
$N_{6,14}^{\alpha,\beta,0}$ &  & no & $\sum_{i=1}^{4}a_{i}d\eta_{i}+a_{5}%
\eta_{3}\wedge\eta_{4}+a_{6}\omega_{1}\wedge\omega_{2}$ & $a_{1}a_{2}a_{5}%
\neq0$\\
$N_{6,15}^{0,\beta,\gamma,0}$ &  & no & $\sum_{i=1}^{4}a_{i}d\eta_{i}%
+a_{1}\gamma^{2}\eta_{2}\wedge\omega_{2}+a_{5}\eta_{3}\wedge\eta_{4}+\ $ &
$a_{1}a_{5}\neq0,\ $\\
&  &  & $+a_{6}\omega_{1}\wedge\omega_{2}$ & $a_{1}a_{6}\neq0$\\
$N_{6,16}^{0,0}$ &  & no & $\sum_{i=1}^{4}a_{i}d\eta_{i}+a_{5}\eta_{3}%
\wedge\eta_{4}+a_{6}\omega_{1}\wedge\omega_{2}$ & $a_{1}a_{5}\neq0$\\
$N_{6,17}^{0}$ &  & no & $\sum_{i=1}^{4}a_{i}d\eta_{i}+a_{5}\eta_{1}\wedge
\eta_{2}+a_{6}\eta_{1}\wedge\omega_{2}+\ $ & $a_{3}a_{5}\neq0,$\\
&  &  & $+a_{7}\eta_{2}\wedge\omega_{2}+a_{8}\omega_{1}\wedge\omega_{2}$ &
$\;a_{4}a_{5}\neq0$\\
$N_{6,18}^{0,\beta,0}$ &  & no & $\sum_{i=1}^{4}a_{i}d\eta_{i}+a_{5}\eta
_{3}\wedge\eta_{4}+a_{6}\omega_{1}\wedge\omega_{2}$ & $a_{2}\left(  a_{1}%
^{2}+a_{2}^{2}\right)  \neq0$\\
$N_{6,20}^{-1,0}$ &  & no & $\sum_{i=1}^{4}a_{i}d\eta_{i}+a_{5}\eta_{2}%
\wedge\eta_{4}+a_{6}\omega_{1}\wedge\eta_{1}+a_{7}\omega_{2}\wedge\eta_{1}$ &
$a_{3}a_{5}a_{6}\neq0$\\
$N_{6,20}^{0,-1}$ &  &  & $\sum_{i=1}^{4}a_{i}d\eta_{i}+a_{5}\eta_{2}%
\wedge\eta_{3}+a_{6}\omega_{1}\wedge\eta_{1}+a_{7}\omega_{2}\wedge\eta_{1}$ &
$a_{4}a_{5}a_{7}\neq0$\\
$N_{6,22}^{\alpha,0}$ &  & no & $\sum_{i=1}^{4}a_{i}d\eta_{i}+a_{5}\eta
_{3}\wedge\eta_{4}+a_{6}\eta_{3}\wedge\omega_{2}+$ & \\
&  &  & $a_{7}\eta_{4}\wedge\omega_{1}+a_{8}\omega_{1}\wedge\omega_{2}$ &
$a_{1}a_{2}a_{6}\neq0$\\
$N_{6,23}^{\alpha,0}$ &  & no & $a_{1}d\eta_{1}+a_{2}d\eta_{2}+a_{3}\eta
_{3}\wedge\omega_{1}+a_{4}\eta_{3}\wedge\omega_{2}+$ & \\
&  &  & $a_{5}\eta_{3}\wedge\eta_{4}+a_{6}\left(  \eta_{1}\wedge\omega
_{1}+\alpha\eta_{4}\wedge\omega_{2}\right)  +a_{7}\omega_{1}\wedge\omega_{2}$%
& $a_{5}\left(  a_{1}^{2}+a_{2}^{2}\right)  \neq0$\\
$N_{6,26}^{\alpha}$ &  & no & $\sum_{i=1}^{4}a_{i}d\eta_{i}+a_{5}\eta
_{3}\wedge\eta_{4}+a_{6}\omega_{1}\wedge\eta_{1}+a_{7}\omega_{2}\wedge\eta
_{1}$ & $a_{2}a_{5}a_{6}\neq0$\\
$N_{6,27}^{0}$ &  & no & $\sum_{i=1}^{4}a_{i}d\eta_{i}+a_{5}\eta_{1}\wedge
\eta_{2}+a_{6}\eta_{1}\wedge\omega_{2}+$ & \\
&  &  & $a_{7}\omega_{1}\wedge\eta_{2}+a_{8}\omega_{1}\wedge\omega_{2}$ &
$a_{5}\left(  a_{3}^{2}+a_{4}^{2}\right)  \neq0$\\
$N_{6,28}$ &  & yes & $\sum_{i=1}^{4}a_{i}d\eta_{i}+a_{5}\omega_{1}%
\wedge\omega_{2}$ & $a_{1}\left(  a_{2}^{2}+a_{3}^{2}\right)  \neq0$\\
$N_{6,29}^{\alpha,\beta}$ & $\alpha\neq0$ & yes & $\sum_{i=1}^{4}a_{i}%
d\eta_{i}+a_{5}\omega_{1}\wedge\omega_{2}+a_{6}\eta_{2}\wedge\eta_{4}$ &
$a_{1}a_{4}\neq0$\\
$N_{6,29}^{0,\beta}$ &  & yes & $\sum_{i=1}^{4}a_{i}d\eta_{i}+a_{5}\omega
_{1}\wedge\omega_{2}+a_{6}\eta_{3}\wedge\eta_{4}$ & $a_{1}a_{4}\neq0$\\
$N_{6,30}^{\alpha}$ &  & yes & $\sum_{i=1}^{4}a_{i}d\eta_{i}+a_{5}\omega
_{1}\wedge\omega_{2}$ & $a_{1}a_{4}\neq0$\\
$N_{6,32}^{\alpha}$ &  & yes & $\sum_{i=1}^{4}a_{i}d\eta_{i}+a_{5}\omega
_{1}\wedge\omega_{2}+a_{6}\omega_{2}\wedge\eta_{4}$ & $a_{1}\neq0$\\
$N_{6,33}$ &  & yes & $\sum_{i=1}^{4}a_{i}d\eta_{i}+a_{5}\omega_{1}%
\wedge\omega_{2}$ & $a_{1}a_{4}\neq0$\\
$N_{6,34}^{\alpha}$ &  & yes & $\sum_{i=1}^{4}a_{i}d\eta_{i}+a_{5}\omega
_{1}\wedge\omega_{2}$ & $a_{1}a_{4}\neq0$\\
$N_{6,35}^{\alpha,\beta}$ & $\alpha\neq0$ & yes & $\sum_{i=1}^{4}a_{i}%
d\eta_{i}+a_{5}\omega_{1}\wedge\omega_{2}$ & $a_{1}a_{4}\neq0$\\
$N_{6,37}^{\alpha}$ &  & yes & $\sum_{i=1}^{4}a_{i}d\eta_{i}+a_{5}\omega
_{1}\wedge\omega_{2}$ & $a_{1}\neq0$\\
$N_{6,38}$ &  & no & $\sum_{i=1}^{4}a_{i}d\eta_{i}+a_{5}\omega_{1}\wedge
\eta_{4}+a_{6}\omega_{2}\wedge\eta_{4}$ & $a_{1}\left(  a_{5}-a_{6}\right)
\neq0$\\
$N_{6,39}$ &  & no & $\sum_{i=1}^{4}a_{i}d\eta_{i}+a_{5}\omega_{1}\wedge
\eta_{4}+a_{6}\omega_{2}\wedge\eta_{4}$ & $a_{1}a_{5}\neq0$\\\hline
\end{tabular}
\end{table}

\subsection*{Acknowledgment}
The author expresses his gratitude to M. Goze and A. Medina for
useful comments. During the preparation of this work, the author
was supported by a research project PR1/05-13283 of the U.C.M..

\section*{Appendix.}

In this appendix we give the Maurer-Cartan equations of the indecomposable solvable non-nilpotent Lie algebras in dimensions three and five, and those of dimension six having a four dimensional nilradical.

The notation and indices for the three dimensional Lie algebras correspond to those given in \cite{Mu2}. In particular, the nilpotent Heisenberg Lie algebra $A_{3,1}$ has been included (see section 3.1). The notation for the five dimensional solvable Lie algebras has been taken from \cite{Mu3}, while the list of six dimensional solvable Lie algebras with four dimensional nilradical has been adapted from \cite{Tu}.

\begin{enumerate}

\item For the three and five dimensional Lie algebras $\left\{\omega_{1},\omega_{2},\omega_{3}\right\}$, respectively $\left\{\omega_{1},..,\omega_{5}\right\}$ denote the dual bases of the algebra. In particular $A_{3,1}$ is nilpotent, and has been included for technical purposes.

\item For the six dimensional Lie algebras with four dimensional nilradical, $\left\{\eta_{1},..,\eta_{4}\right\}$ denotes the dual basis of the nilradical, while $\left\{\omega_{1},\omega_{2}\right\}$ is a dual basis of the space of nil-independent elements (i.e., linearly independent non-nilpotent derivations of the nilradical). For all these Lie algebras $d\omega_{i}=0$.

\item The restrictions on the parameters of the Lie algebras have been indicated after the Maurer-Cartan equations. We remark that some authors alter the numbering of the isomorphism classes in references \cite{Mu2,Mu3} for special values of the parameters.

\end{enumerate}

\begin{table}[h]
\caption{Indecomposable solvable Lie algebras of dimension 3}
\begin{tabular}
[c]{lllll}%
Name & \multicolumn{4}{c}{Maurer-Cartan equations}\\\hline
$A_{3,1}$ & $d\omega_{1}=\omega_{2}\wedge\omega_{3},$ & $d\omega_{2}%
=d\omega_{3}=0,$ & (nilpotent) & \\
$A_{3,2}$ & $d\omega_{1}=\omega_{1}\wedge\omega_{3}+\omega_{2}\wedge\omega
_{3},$ & $d\omega_{2}=\omega_{2}\wedge\omega_{3},$ & $d\omega_{3}=0$ & \\
$A_{3,3}$ & $d\omega_{1}=\omega_{1}\wedge\omega_{3},$ & $d\omega_{2}%
=\omega_{2}\wedge\omega_{3},$ & $d\omega_{3}=0$ & \\
$A_{3,4}^{\alpha}$ & $d\omega_{1}=\omega_{1}\wedge\omega_{3},$ & $d\omega
_{2}=\alpha\omega_{2}\wedge\omega_{3},$ & $d\omega_{3}=0,$ & $-1\leq\alpha
\leq1,\alpha\neq0$\\
$A_{3,5}^{p}$ & $d\omega_{1}=p\omega_{1}\wedge\omega_{3}+\omega_{2}%
\wedge\omega_{3},$ & $d\omega_{2}=-\omega_{1}\wedge\omega_{3}+p\omega
_{2}\wedge\omega_{3}$ & $d\omega_{3}=0,$ & $p\geq0$\\ \hline
\end{tabular}
\end{table}

\begin{table}
\caption{Indecomposable solvable Lie algebras of dimension 5}
\begin{tabular}
[c]{lllll}%
Name & \multicolumn{4}{c}{Maurer-Cartan equations}\\\hline
$\frak{g}_{5,7}^{\alpha\beta\gamma}$ & $d\omega_{1}=\omega_{1}\wedge\omega
_{5},$ & $d\omega_{2}=\alpha\omega_{2}\wedge\omega_{5},$ & $d\omega_{3}%
=\beta\omega_{3}\wedge\omega_{5},$ & $d\omega_{4}=\gamma\omega_{4}\wedge
\omega_{5}$\\
& $d\omega_{5}=0,$ & \multicolumn{2}{l}{$\left(  -1\leq\alpha,\beta,\gamma
\leq1,\;\alpha\gamma\beta\neq0\right)  $} & \\
$\frak{g}_{5,8}^{\gamma}$ & $d\omega_{1}=\omega_{2}\wedge\omega_{5},$ &
$d\omega_{2}=0,$ & $d\omega_{3}=\omega_{3}\wedge\omega_{5}$ & $d\omega
_{4}=\gamma\omega_{4}\wedge\omega_{5}$\\
& $d\omega_{5}=0$ & $\left(  0<\left|  \gamma\right|  \leq1\right)  $ &  & \\
$\frak{g}_{5,9}^{\gamma\beta}$ & \multicolumn{2}{l}{$d\omega_{1}=\omega
_{1}\wedge\omega_{5}+\omega_{2}\wedge\omega_{5},$} & $d\omega_{2}=\omega
_{2}\wedge\omega_{5},$ & $d\omega_{3}=\beta\omega_{3}\wedge\omega_{5},$\\
& $d\omega_{4}=\gamma\omega_{4}\wedge\omega_{5},$ & $d\omega_{5}=0$ & $\left(
0\neq\gamma\leq\beta\right)  $ & \\
$\frak{g}_{5,10}$ & $d\omega_{1}=\omega_{2}\wedge\omega_{5},$ & $d\omega
_{2}=\omega_{3}\wedge\omega_{5},$ & $d\omega_{3}=0,$ & $d\omega_{4}=\omega
_{4}\wedge\omega_{5},$\\
& $d\omega_{5}=0.$ &  &  & \\
$\frak{g}_{5,11}^{\gamma}$ & \multicolumn{2}{l}{$d\omega_{1}=\omega_{1}%
\wedge\omega_{5}+\omega_{2}\wedge\omega_{5},$} & \multicolumn{2}{l}{$d\omega
_{2}=\omega_{2}\wedge\omega_{5}+\omega_{3}\wedge\omega_{5},$}\\
& $d\omega_{3}=\omega_{3}\wedge\omega_{5},$ & $d\omega_{4}=\gamma\omega
_{4}\wedge\omega_{5},$ & $d\omega_{5}=0$ & $\left(  \gamma\neq0\right)  $\\
$\frak{g}_{5,12}$ & \multicolumn{2}{l}{$d\omega_{1}=\omega_{1}\wedge\omega
_{5}+\omega_{2}\wedge\omega_{5},$} & \multicolumn{2}{l}{$d\omega_{2}%
=\omega_{2}\wedge\omega_{5}+\omega_{3}\wedge\omega_{5},$}\\
& \multicolumn{2}{l}{$d\omega_{3}=\omega_{3}\wedge\omega_{5}+\omega_{4}%
\wedge\omega_{5},$} & $d\omega_{4}=\omega_{4}\wedge\omega_{5},$ & $d\omega
_{5}=0.$\\
$\frak{g}_{5,13}^{\gamma,p,s}$ & $d\omega_{1}=\omega_{1}\wedge\omega_{5},$ &
$d\omega_{2}=\gamma\omega_{2}\wedge\omega_{5},$ & \multicolumn{2}{l}{$d\omega
_{3}=p\omega_{3}\wedge\omega_{5}+s\omega_{4}\wedge\omega_{5},$}\\
& \multicolumn{2}{l}{$d\omega_{4}=-s\omega_{3}\wedge\omega_{5}+p\omega
_{4}\wedge\omega_{5},$} & $d\omega_{5}=0$ & $\left(  \left|  \gamma\right|
\leq1,\;\gamma s\neq0\right)  $\\
$\frak{g}_{5,14}^{p}$ & $d\omega_{1}=\ \omega_{2}\wedge\omega_{5},$ &
$d\omega_{2}=0,$ & \multicolumn{2}{l}{$d\omega_{3}=p\omega_{3}\wedge\omega
_{5}+\omega_{4}\wedge\omega_{5},$}\\
& \multicolumn{2}{l}{$d\omega_{4}=-\omega_{3}\wedge\omega_{5}+p\omega
_{4}\wedge\omega_{5},$} & $d\omega_{5}=0.$ & \\
$\frak{g}_{5,15}^{\gamma}$ & \multicolumn{2}{l}{$d\omega_{1}=\omega_{1}%
\wedge\omega_{5}+\omega_{2}\wedge\omega_{5},$} & $d\omega_{2}=\omega_{2}%
\wedge\omega_{5},$ & $d\omega_{3}=\gamma\omega_{3}\wedge\omega_{5}+\omega
_{4}\wedge\omega_{5},$\\
& $d\omega_{4}=\gamma\omega_{4}\wedge\omega_{5},$ & $d\omega_{5}=0$ & $\left(
-1\leq\gamma\leq1\right)  $ & \\
$\frak{g}_{5,16}^{p,s}$ & \multicolumn{2}{l}{$d\omega_{1}=\omega_{1}%
\wedge\omega_{5}+\omega_{2}\wedge\omega_{5},$} & $d\omega_{1}=\omega_{2}%
\wedge\omega_{5},$ & $d\omega_{3}=p\omega_{3}\wedge\omega_{5}+s\omega
_{4}\wedge\omega_{5},$\\
& \multicolumn{2}{l}{$d\omega_{4}=-s\omega_{3}\wedge\omega_{5}+p\omega
_{4}\wedge\omega_{5},$} & $d\omega_{5}=0$ & $\left(  s\neq0\right)  $\\
$\frak{g}_{5,17}^{pqs}$ & \multicolumn{2}{l}{$d\omega_{1}=p\omega_{1}%
\wedge\omega_{5}+\omega_{2}\wedge\omega_{5},$} & \multicolumn{2}{l}{$d\omega
_{2}=-\omega_{1}\wedge\omega_{5}+p\omega_{2}\wedge\omega_{5},$}\\
& \multicolumn{2}{l}{$d\omega_{3}=q\omega_{3}\wedge\omega_{5}+s\omega
_{4}\wedge\omega_{5},$} & \multicolumn{2}{l}{$d\omega_{4}=-s\omega_{3}%
\wedge\omega_{5}+q\omega_{4}\wedge\omega_{5},$}\\
& $d\omega_{5}=0$ & $\left(  s\neq0\right)  $ &  & \\
$\frak{g}_{5,18}^{p}$ & \multicolumn{2}{l}{$d\omega_{1}=p\omega_{1}%
\wedge\omega_{5}+\omega_{2}\wedge\omega_{5}+\omega_{3}\wedge\omega_{5},$} &
\multicolumn{2}{l}{$d\omega_{2}=-\omega_{1}\wedge\omega_{5}+p\omega_{2}%
\wedge\omega_{5}+\omega_{4}\wedge\omega_{5},$}\\
& \multicolumn{2}{l}{$d\omega_{3}=p\omega_{3}\wedge\omega_{5}+\omega_{4}%
\wedge\omega_{5},$} & \multicolumn{2}{l}{$d\omega_{4}=-\omega_{3}\wedge
\omega_{5}-p\omega_{4}\wedge\omega_{5},$}\\
& $d\omega_{5}=0$ & $\left(  p\geq0\right)  $ &  & \\
$\frak{g}_{5,19}^{\alpha\beta}$ & \multicolumn{2}{l}{$d\omega_{1}=\omega
_{2}\wedge\omega_{3}+\left(  1+\alpha\right)  \omega_{1}\wedge\omega_{5},$} &
$d\omega_{2}=\omega_{2}\wedge\omega_{5},$ & $d\omega_{3}=\alpha\omega
_{3}\wedge\omega_{5},$\\
& $d\omega_{4}=\beta\omega_{4}\wedge\omega_{5},$ & $\left(  \beta\neq0\right)
$ &  & \\
$\frak{g}_{5,20}^{\alpha}$ & \multicolumn{2}{l}{$d\omega_{1}=\omega_{2}%
\wedge\omega_{3}+\left(  1+\alpha\right)  \omega_{1}\wedge\omega_{5}%
+\omega_{4}\wedge\omega_{5},$} & $d\omega_{2}=\omega_{2}\wedge\omega_{5},$ &
$d\omega_{3}=\alpha\omega_{3}\wedge\omega_{5},$\\
& $d\omega_{4}=\left(  1+\alpha\right)  \omega_{4}\wedge\omega_{5},$ &
$d\omega_{5}=0$ &  & \\
$\frak{g}_{5,21}$ & \multicolumn{2}{l}{$d\omega_{1}=\omega_{2}\wedge\omega
_{3}+2\omega_{1}\wedge\omega_{5},$} & \multicolumn{2}{l}{$d\omega_{2}%
=\omega_{2}\wedge\omega_{5},$}\\
& \multicolumn{2}{l}{$d\omega_{3}=\omega_{2}\wedge\omega_{5}+\omega_{3}%
\wedge\omega_{5},$} & \multicolumn{2}{l}{$d\omega_{4}=\omega_{3}\wedge
\omega_{5}+\omega_{4}\wedge\omega_{5},$}\\
& $d\omega_{5}=0.$ &  &  & \\
$\frak{g}_{5,22}$ & $d\omega_{1}=\omega_{2}\wedge\omega_{3},$ & $d\omega
_{2}=0,$ & $d\omega_{3}=\omega_{2}\wedge\omega_{5},$ & $d\omega_{4}=\omega
_{4}\wedge\omega_{5},$\\
& $d\omega_{5}=0.$ &  &  & \\
$\frak{g}_{5,23}^{\beta}$ & $d\omega_{1}=\omega_{2}\wedge\omega_{3}%
+2\omega_{1}\wedge\omega_{5},$ & $d\omega_{2}=\omega_{2}\wedge\omega_{5},$ &
\multicolumn{2}{l}{$d\omega_{3}=\omega_{3}\wedge\omega_{5}+\omega_{2}%
\wedge\omega_{5},$}\\
& $d\omega_{4}=\beta\omega_{4}\wedge\omega_{5},$ & $d\omega_{5}=0$ & $\left(
\beta\neq0\right)  $ & \\
$\frak{g}_{5,24}$ & \multicolumn{2}{l}{$d\omega_{1}=\omega_{2}\wedge\omega
_{3}+2\omega_{1}\wedge\omega_{5}+\varepsilon\omega_{4}\wedge\omega_{5},$} &
$d\omega_{2}=\omega_{2}\wedge\omega_{5},$ & $\left(  \varepsilon=\pm1\right)
$\\
& \multicolumn{2}{l}{$d\omega_{3}=\omega_{2}\wedge\omega_{5}+\omega_{3}%
\wedge\omega_{5},$} & $d\omega_{4}=2\omega_{4}\wedge\omega_{5},$ &
$d\omega_{5}=0.$\\
$\frak{g}_{5,25}^{\beta,p}$ & \multicolumn{2}{l}{$d\omega_{1}=\omega_{2}%
\wedge\omega_{3}+2p\omega_{1}\wedge\omega_{5},$} & \multicolumn{2}{l}{$d\omega
_{2}=p\omega_{2}\wedge\omega_{5}-\omega_{3}\wedge\omega_{5},$}\\
& \multicolumn{2}{l}{$d\omega_{3}=\omega_{2}\wedge\omega_{5}+p\omega_{3}%
\wedge\omega_{5},$} & $d\omega_{4}=\beta\omega_{4}\wedge\omega_{5},$ &
$d\omega_{5}=0\;\left(  \beta\neq0\right)  $\\\hline
\end{tabular}
\end{table}

\bigskip%

\begin{table}
\caption{Indecomposable solvable Lie algebras of dimension 5 (cont.)}
\begin{tabular}
[c]{lllll}%
Name & \multicolumn{4}{c}{Maurer-Cartan equations}\\\hline
$\frak{g}_{5,26}^{p\varepsilon}$ & \multicolumn{2}{l}{$d\omega_{1}=\omega
_{2}\wedge\omega_{3}+2p\omega_{1}\wedge\omega_{5}+\varepsilon\omega_{4}%
\wedge\omega_{5},$} & \multicolumn{2}{l}{$d\omega_{2}=p\omega_{2}\wedge
\omega_{5}-\omega_{3}\wedge\omega_{5},$}\\
& \multicolumn{2}{l}{$d\omega_{3}=\omega_{2}\wedge\omega_{5}+p\omega_{3}%
\wedge\omega_{5},$} & $d\omega_{4}=2p\omega_{4}\wedge\omega_{5},$ &
$d\omega_{5}=0\;\left(  \varepsilon=\pm1\right)  $\\
$\frak{g}_{5,27}$ & \multicolumn{2}{l}{$d\omega_{1}=\omega_{2}\wedge\omega
_{3}+\omega_{1}\wedge\omega_{5}+\omega_{4}\wedge\omega_{5},$} & $d\omega
_{2}=0,$ & $d\omega_{3}=\omega_{3}\wedge\omega_{5},$\\
& \multicolumn{2}{l}{$d\omega_{4}=\omega_{3}\wedge\omega_{5}+\omega_{4}%
\wedge\omega_{5},$} & $d\omega_{5}=0.$ & \\
$\frak{g}_{5,28}^{\alpha}$ & \multicolumn{2}{l}{$d\omega_{1}=\omega_{2}%
\wedge\omega_{3}+\left(  1+\alpha\right)  \omega_{1}\wedge\omega_{5},$} &
$d\omega_{2}=\alpha\omega_{2}\wedge\omega_{5},$ & $d\omega_{3}=\omega
_{3}\wedge\omega_{5},$\\
& \multicolumn{2}{l}{$d\omega_{4}=\omega_{3}\wedge\omega_{5}+\omega_{4}%
\wedge\omega_{5},$} & $d\omega_{5}=0.$ & \\
$\frak{g}_{5,29}$ & \multicolumn{2}{l}{$d\omega_{1}=\omega_{2}\wedge\omega
_{3}+\omega_{1}\wedge\omega_{5},$} & $d\omega_{2}=\omega_{2}\wedge\omega_{5},$%
& $d\omega_{3}=0,$\\
& $d\omega_{4}=\omega_{3}\wedge\omega_{5},$ & $d\omega_{5}=0.$ &  & \\
$\frak{g}_{5,30}^{\alpha}$ & \multicolumn{2}{l}{$d\omega_{1}=\omega_{2}%
\wedge\omega_{4}+\left(  2+\alpha\right)  \omega_{1}\wedge\omega_{5},$} &
\multicolumn{2}{l}{$d\omega_{2}=\omega_{3}\wedge\omega_{4}+\left(
1+\alpha\right)  \omega_{2}\wedge\omega_{5},$}\\
& $d\omega_{3}=\omega\alpha_{3}\wedge\omega_{5},$ & $d\omega_{4}=\omega
_{4}\wedge\omega_{5},$ & $d\omega_{5}=0.$ & \\
$\frak{g}_{5,31}$ & \multicolumn{2}{l}{$d\omega_{1}=\omega_{2}\wedge\omega
_{4}+3\omega_{1}\wedge\omega_{5},$} & \multicolumn{2}{l}{$d\omega_{1}%
=\omega_{3}\wedge\omega_{4}+2\omega_{2}\wedge\omega_{5},$}\\
& \multicolumn{2}{l}{$d\omega_{3}=\omega_{3}\wedge\omega_{5}+\omega_{4}%
\wedge\omega_{5},$} & $d\omega_{4}=\omega_{4}\wedge\omega_{5},$ & $d\omega
_{5}=0.$\\
$\frak{g}_{5,32}^{\alpha}$ & \multicolumn{2}{l}{$d\omega_{1}=\omega_{2}%
\wedge\omega_{4}+\omega_{1}\wedge\omega_{5}+\alpha\omega_{3}\wedge\omega_{5}%
,$} & \multicolumn{2}{l}{$d\omega_{2}=\omega_{3}\wedge\omega_{4}+\omega
_{2}\wedge\omega_{5},$}\\
& $d\omega_{3}=\omega_{3}\wedge\omega_{5},$ & $d\omega_{4}=0,$ & $d\omega
_{5}=0.$ & \\
$\frak{g}_{5,33}^{\beta\gamma}$ & $d\omega_{1}=\omega_{1}\wedge\omega_{4},$ &
$d\omega_{2}=\omega_{2}\wedge\omega_{5},$ & \multicolumn{2}{l}{$d\omega
_{3}=\beta\omega_{3}\wedge\omega_{4}+\gamma\omega_{3}\wedge\omega_{5},$}\\
& $d\omega_{4}=0,$ & $d\omega_{5}=0$ & $\left(  \beta^{2}+\gamma^{2}\right)
\neq0.$ & \\
$\frak{g}_{5,34}^{\alpha}$ & \multicolumn{2}{l}{$d\omega_{1}=\alpha\omega
_{1}\wedge\omega_{4}+\omega_{1}\wedge\omega_{5},$} &
\multicolumn{2}{l}{$d\omega_{2}=\omega_{2}\wedge\omega_{4}+\omega_{3}%
\wedge\omega_{5},$}\\
& $d\omega_{3}=\omega_{3}\wedge\omega_{4},$ & $d\omega_{4}=0,$ & $d\omega
_{5}=0.$ & \\
$\frak{g}_{5,35}^{\alpha,\beta}$ & \multicolumn{2}{l}{$d\omega_{1}=\beta
\omega_{1}\wedge\omega_{4}+\alpha\omega_{1}\wedge\omega_{5},$} &
\multicolumn{2}{l}{$d\omega_{2}=\omega_{2}\wedge\omega_{4}+\omega_{3}%
\wedge\omega_{5},$}\\
& \multicolumn{2}{l}{$d\omega_{3}=\omega_{3}\wedge\omega_{4}-\omega_{2}%
\wedge\omega_{5},$} & $d\omega_{4}=d\omega_{5}=0,$ & $\left(  \alpha^{2}%
+\beta^{2}\neq0\right)  $\\
$\frak{g}_{5,36}$ & \multicolumn{2}{l}{$d\omega_{1}=\omega_{2}\wedge\omega
_{3}+\omega_{1}\wedge\omega_{4},$} & \multicolumn{2}{l}{$d\omega_{2}%
=\omega_{2}\wedge\omega_{4}-2\omega_{2}\wedge\omega_{5},$}\\
& $d\omega_{3}=\omega_{3}\wedge\omega_{5},$ & $d\omega_{4}=d\omega_{5}=0.$ &
& \\
$\frak{g}_{5,37}$ & \multicolumn{2}{l}{$d\omega_{1}=\omega_{2}\wedge\omega
_{3}+2\omega_{1}\wedge\omega_{4},$} & \multicolumn{2}{l}{$d\omega_{2}%
=\omega_{2}\wedge\omega_{4}+\omega_{3}\wedge\omega_{5},$}\\
& \multicolumn{2}{l}{$d\omega_{3}=\omega_{3}\wedge\omega_{4}-\omega_{2}%
\wedge\omega_{5},$} & $d\omega_{4}=d\omega_{5}=0.$ & \\
$\frak{g}_{5,38}$ & $d\omega_{1}=\omega_{1}\wedge\omega_{4},$ & $d\omega
_{2}=\omega_{2}\wedge\omega_{5},$ & $d\omega_{3}=\omega_{4}\wedge\omega_{5},$%
& $d\omega_{4}=d\omega_{5}=0.$\\
$\frak{g}_{5,39}$ & \multicolumn{2}{l}{$d\omega_{1}=\omega_{1}\wedge\omega
_{4}+\omega_{2}\wedge\omega_{5},$} & \multicolumn{2}{l}{$d\omega_{2}%
=\omega_{2}\wedge\omega_{4}-\omega_{1}\wedge\omega_{5},$}\\
& $d\omega_{3}=\omega_{4}\wedge\omega_{5},$ & $d\omega_{4}=d\omega_{5}=0.$ &
&\\\hline
\end{tabular}
\end{table}

\begin{table}
\caption{Indecomposable solvable Lie algebras of dimension 6 with $\dim(NR)=4$.}
\begin{tabular}
[c]{lllll}%
Name & \multicolumn{4}{l}{Maurer-Cartan equations}\\\hline
$N_{6,1}^{\alpha\beta\gamma\delta}$ & $d\eta_{1}=\alpha\omega_{1}\wedge
\eta_{1}+\beta\omega_{2}\wedge\eta_{1},$ & $d\eta_{2}=\gamma\omega_{1}%
\wedge\eta_{2}+\delta\omega_{2}\wedge\eta_{2},$ & $\gamma^{2}+\delta^{2}\neq0$%
& \\
$\left(  \alpha\beta\neq0\right)  $ & $d\eta_{3}=\omega_{2}\wedge\eta_{3},$ &
$d\eta_{4}=\omega_{1}\wedge\eta_{4},$ & $\ d\omega_{1}=d\omega_{2}=0$ & \\
$N_{6,2}^{\alpha\beta\gamma}$ & $d\eta_{1}=\alpha\omega_{1}\wedge\eta
_{1}+\beta\omega_{2}\wedge\eta_{1},\ $ & $d\eta_{2}=\omega_{1}\wedge\eta
_{2}+\gamma\omega_{2}\wedge\eta_{2},\ $ & $\ $ & \\
$\alpha^{2}+\beta^{2}\neq0$ & $d\eta_{3}=\omega_{2}\wedge\eta_{3},\ $ &
$d\eta_{4}=\omega_{1}\wedge\eta_{3}+\omega_{2}\wedge\eta_{4},\ $ &
$d\omega_{1}=d\omega_{2}=0\ $ & \\
$N_{6,3}^{\alpha}$ & $d\eta_{1}=\omega_{1}\wedge\eta_{1}+\alpha\omega
_{2}\wedge\eta_{1},\ $ & $d\eta_{2}=\omega_{2}\wedge\eta_{1}+\omega_{1}%
\wedge\eta_{2}+\alpha\omega_{2}\wedge\eta_{2},\ $ & $\ $ & \\
& $d\eta_{3}=\omega_{2}\wedge\eta_{3},\ $ & $d\eta_{4}=\omega_{2}\wedge
\eta_{4}+\omega_{1}\wedge\eta_{3},\ $ & $d\omega_{1}=d\omega_{2}=0\ $ & \\
$N_{6,4}^{\alpha\beta}$ & $d\eta_{1}=\omega_{1}\wedge\eta_{1}-\omega_{2}%
\wedge\eta_{2},\ $ & $d\eta_{2}=\omega_{2}\wedge\eta_{1}+\omega_{1}\wedge
\eta_{2},$ &  & \\
$\alpha\neq0\ $ & $d\eta_{3}=\alpha\omega_{2}\wedge\eta_{3},\ $ & $d\eta
_{4}=\beta\omega_{2}\wedge\eta_{3}+\alpha\omega_{2}\wedge\eta_{4}+\omega
_{1}\wedge\eta_{3},\ $ & $d\omega_{1}=d\omega_{2}=0\ $ & \\
$N_{6,5}^{\alpha\beta}$ & $d\eta_{1}=\alpha\omega_{1}\wedge\eta_{1}%
+\beta\omega_{2}\wedge\eta_{1},\ $ & $d\eta_{2}=\omega_{2}\wedge\eta_{2},\ $ &
$\ $ & \\
$\alpha\beta\neq0$ & $d\eta_{3}=\omega_{1}\wedge\eta_{3},\ $ & $d\eta
_{4}=\omega_{1}\wedge\eta_{3}+\omega_{1}\wedge\eta_{4},\ $ & $d\omega
_{1}=d\omega_{2}=0$ & \\
$N_{6,6}^{\alpha\beta}$ & $\ d\eta_{1}=\alpha\omega_{1}\wedge\eta_{1}%
+\omega_{2}\wedge\eta_{1},\ $ & $\ d\eta_{2}=\omega_{2}\wedge\eta_{1}%
+\alpha\omega_{1}\wedge\eta_{2}+\omega_{2}\wedge\eta_{2},\ $ & $\ $ & \\
$\alpha^{2}+\beta^{2}\neq0$ & $\ d\eta_{3}=\omega_{1}\wedge\eta_{3},$ &
$\ d\eta_{4}=\omega_{1}\wedge\eta_{3}+\omega_{1}\wedge\eta_{4}-\beta\omega
_{2}\wedge\eta_{3},$ & $d\omega_{1}=d\omega_{2}=0\ $ & \\
$N_{6,7}^{\alpha\beta\gamma}$ & $d\eta_{1}=\alpha\omega_{1}\wedge\eta
_{1}-\omega_{2}\wedge\eta_{2}+\gamma\omega_{2}\wedge\eta_{1},\ $ & $d\eta
_{2}=\alpha\omega_{1}\wedge\eta_{2}+\omega_{2}\wedge\eta_{1}+\gamma\omega
_{2}\wedge\eta_{2},\ $ &  & \\
$\alpha^{2}+\beta^{2}\neq0$ & $d\eta_{3}=\omega_{1}\wedge\eta_{3},\ $ &
$d\eta_{4}=\omega_{1}\wedge\eta_{3}+\omega_{1}\wedge\eta_{4}+\beta\omega
_{2}\wedge\eta_{3},\ $ & $d\omega_{1}=d\omega_{2}=0\ $ & \\
$N_{6,8}$ & $d\eta_{1}=\omega_{1}\wedge\eta_{1},$ & $d\eta_{2}=\omega
_{2}\wedge\eta_{2},$ &  & \\
& $d\eta_{3}=\omega_{2}\wedge\eta_{3},$ & $d\eta_{4}=\omega_{2}\wedge\eta
_{3}+\omega_{2}\wedge\eta_{4},$ & $d\omega_{1}=d\omega_{2}=0$ & \\
$N_{6,9}^{\alpha}$ & $d\eta_{4}=\omega_{1}\wedge\eta_{1},$ & $d\eta_{2}%
=\omega_{2}\wedge\eta_{2},$ &  & \\
& $d\eta_{3}=\omega_{2}\wedge\eta_{2}+\omega_{2}\wedge\eta_{3},$ & $d\eta
_{4}=\omega_{1}\wedge\eta_{2}+\alpha\omega_{2}\wedge\eta_{3}+\omega_{2}%
\wedge\eta_{4},$ & $d\omega_{1}=d\omega_{2}=0$ & \\
$N_{6,10}^{\alpha\beta}$ & $d\eta_{1}=\alpha\omega_{1}\wedge\eta_{1}%
+\omega_{1}\wedge\eta_{2}+\omega_{2}\wedge\eta_{1},$ & $d\eta_{2}=0,$ &  & \\
& $d\eta_{3}=\omega_{1}\wedge\eta_{3}+\omega_{2}\wedge\eta_{2},$ & $d\eta
_{4}=\beta\omega_{1}\wedge\eta_{2}+\omega_{1}\wedge\eta_{4}+\omega_{2}%
\wedge\eta_{3},$ & $d\omega_{1}=d\omega_{2}=0$ & \\
$N_{6,11}^{\alpha}$ & $d\eta_{1}=\omega_{2}\wedge\eta_{1},$ & $d\eta
_{2}=\omega_{1}\wedge\eta_{1}+\omega_{2}\wedge\eta_{2},$ &  & \\
& $d\eta_{3}=\omega_{1}\wedge\eta_{3}+\alpha\omega_{2}\wedge\eta_{3},$ &
$d\eta_{4}=\omega_{1}\wedge\eta_{3}+\omega_{1}\wedge\eta_{4}+\alpha\omega
_{2}\wedge\eta_{4},$ & $d\omega_{1}=d\omega_{2}=0$ & \\
$N_{6,12}^{\alpha\beta}$ & $d\eta_{1}=\omega_{1}\wedge\eta_{1}-\omega
_{2}\wedge\eta_{3},$ & $d\eta_{3}=\omega_{1}\wedge\eta_{3}+\omega_{2}%
\wedge\eta_{1},$ &  & \\
& \multicolumn{2}{l}{$d\eta_{2}=\omega_{1}\wedge\eta_{1}+\omega_{1}\wedge
\eta_{2}+\alpha\omega_{2}\wedge\eta_{1}+\beta\omega_{2}\wedge\eta_{3}%
-\omega_{2}\wedge\eta_{4},$} &  & \\
& \multicolumn{2}{l}{$d\eta_{4}=\omega_{1}\wedge\eta_{3}+\omega_{1}\wedge
\eta_{4}-\beta\omega_{2}\wedge\eta_{1}+\omega_{2}\wedge\eta_{2}+\alpha
\omega_{2}\wedge\eta_{3},$} & $d\omega_{1}=d\omega_{2}=0$ & \\
$N_{6,13}^{\alpha\beta\gamma\delta}$ & $d\eta_{1}=\alpha\omega_{1}\wedge
\eta_{1}+\beta\omega_{2}\wedge\eta_{1},$ & $d\eta_{2}=\gamma\omega_{1}%
\wedge\eta_{2}+\gamma\omega_{2}\wedge\eta_{2},$ & $\gamma^{2}+\delta^{2}\neq0$%
& \\
$\alpha^{2}+\beta^{2}\neq0$ & $d\eta_{3}=-\omega_{1}\wedge\eta_{4}+\omega
_{2}\wedge\eta_{3},$ & $d\eta_{4}=\omega_{1}\wedge\eta_{3}+\omega_{2}%
\wedge\eta_{4},$ & $d\omega_{1}=d\omega_{2}=0$ & \\
$N_{6,14}^{\alpha\beta\gamma}$ & $d\eta_{1}=\alpha\omega_{1}\wedge\eta
_{1}+\beta\omega_{2}\wedge\eta_{1},$ & $d\eta_{2}=\omega_{2}\wedge\eta_{2},$ &
& \\
$\alpha\beta\neq0$ & $d\eta_{3}=\gamma\omega_{1}\wedge\eta_{3}-\omega
_{1}\wedge\eta_{4},$ & $d\eta_{4}=\omega_{1}\wedge\eta_{3}+\gamma\omega
_{1}\wedge\eta_{4},$ & $d\omega_{1}=d\omega_{2}=0$ & \\
$N_{6,15}^{\alpha\beta\gamma\delta}$ & $d\eta_{1}=\omega_{1}\wedge\eta
_{1}+\gamma\omega_{2}\wedge\eta_{1}-\omega_{2}\wedge\eta_{2},$ & $d\eta
_{2}=\omega_{1}\wedge\eta_{2}+\omega_{2}\wedge\eta_{1}+\gamma\omega_{2}%
\wedge\eta_{2},$ &  & \\
$\beta\neq0$ & $d\eta_{3}=\alpha\omega_{1}\wedge\eta_{3}-\beta\omega_{1}%
\wedge\eta_{4}+\delta\omega_{2}\wedge\eta_{3},$ & $d\eta_{4}=\beta\omega
_{1}\wedge\eta_{3}+\alpha\omega_{1}\wedge\eta_{4}+\delta\omega_{2}\wedge
\eta_{4},$ & $d\omega_{1}=d\omega_{2}=0$ & \\
$N_{6,16}^{\alpha\beta}$ & $d\eta_{4}=\omega_{2}\wedge\eta_{1},$ & $d\eta
_{2}=\omega_{1}\wedge\eta_{1}+\omega_{2}\wedge\eta_{2},$ &  & \\
& $d\eta_{3}=\alpha\omega_{1}\wedge\eta_{3}-\omega_{1}\wedge\eta_{4}%
+\beta\omega_{2}\wedge\eta_{3},$ & $d\eta_{4}=\omega_{1}\wedge\eta_{3}%
+\alpha\omega_{1}\wedge\eta_{4}+\beta\omega_{2}\wedge\eta_{4},$ & $d\omega
_{1}=d\omega_{2}=0$ & \\
$N_{6,17}^{\alpha}$ & $d\eta_{1}=\alpha\omega_{1}\wedge\eta_{1},$ & $d\eta
_{2}=\omega_{1}\wedge\eta_{1}+\alpha\omega_{1}\wedge\eta_{2},$ &  & \\
& $d\eta_{3}=-\omega_{1}\wedge\eta_{4}+\omega_{2}\wedge\eta_{3},$ & $d\eta
_{4}=\omega_{1}\wedge\eta_{3}+\omega_{2}\wedge\eta_{4},$ & $d\omega
_{1}=d\omega_{2}=0$ & \\
$N_{6,18}^{\alpha\beta\gamma}$ & $d\eta_{1}=-\omega_{1}\wedge\eta_{2}%
+\omega_{2}\wedge\eta_{1},$ & $d\eta_{2}=\omega_{1}\wedge\eta_{1}+\omega
_{2}\wedge\eta_{2},$ &  & \\
$\beta\neq0$ & $d\eta_{3}=\alpha\omega_{1}\wedge\eta_{3}-\beta\omega_{1}%
\wedge\eta_{4}+\gamma\omega_{2}\wedge\eta_{3},$ & $d\eta_{4}=\beta\omega
_{1}\wedge\eta_{3}+\alpha\omega_{1}\wedge\eta_{4}+\gamma\omega_{2}\wedge
\eta_{4},$ & $d\omega_{1}=d\omega_{2}=0$ & \\
$N_{6,19}$ & $d\eta_{1}=-\omega_{1}\wedge\eta_{2}+\omega_{2}\wedge\eta_{1},$ &
$d\eta_{2}=\omega_{1}\wedge\eta_{1}+\omega_{2}\wedge\eta_{2},$ &  & \\
& $d\eta_{3}=\omega_{1}\wedge\eta_{1}-\omega_{1}\wedge\eta_{4}+\omega
_{2}\wedge\eta_{3},$ & $d\eta_{4}=\omega_{1}\wedge\eta_{2}+\omega_{1}%
\wedge\eta_{3}+\omega_{2}\wedge\eta_{4},$ & $d\omega_{1}=d\omega_{2}=0$ &\\\hline
\end{tabular}
\end{table}

\bigskip%

\begin{table*}
\caption{Indecomposable solvable Lie algebras of dimension 6 with $\dim(NR)=4$ (cont.)}
\begin{tabular}
[c]{lllll}%
Name & \multicolumn{4}{l}{Maurer-Cartan equations}\\\hline
$N_{6,20}^{\alpha\beta}$ & $d\eta_{1}=\omega_{1}\wedge\omega_{2},$ &
$d\eta_{2}=\alpha\omega_{1}\wedge\eta_{2}+\beta\omega_{2}\wedge\eta_{2},$ &  &
\\
$\alpha^{2}+\beta^{2}\neq0$ & $d\eta_{3}=\omega_{2}\wedge\eta_{3},$ &
$d\eta_{4}=\omega_{1}\wedge\eta_{4},$ & $d\omega_{1}=d\omega_{2}=0$ & \\
$N_{6,21}^{\alpha}$ & $d\eta_{1}=\omega_{1}\wedge\omega_{2},$ & $d\eta
_{2}=\omega_{1}\wedge\eta_{2}+a\omega_{2}\wedge\eta_{2},$ & $\ $ & \\
& $d\eta_{3}=\omega_{2}\wedge\eta_{3},$ & $d\eta_{4}=\omega_{1}\wedge\eta
_{3}+\omega_{2}\wedge\eta_{4},$ & $d\omega_{1}=d\omega_{2}=0\ $ & \\
$N_{6,22}^{\alpha\varepsilon}$ & $d\eta_{1}=\omega_{1}\wedge\eta_{1}%
+\alpha\omega_{2}\wedge\eta_{1},$ & $d\eta_{2}=\omega_{2}\wedge\eta_{2},$ &
$\varepsilon=0,1$ & \\
$\alpha^{2}+\varepsilon^{2}\neq0$ & $d\eta_{3}=\varepsilon\omega_{1}%
\wedge\omega_{2},$ & $d\eta_{4}=\omega_{1}\wedge\eta_{3},$ & $d\omega
_{1}=d\omega_{2}=0\ $ & \\
$N_{6,23}^{\alpha\varepsilon}$ & $d\eta_{1}=\omega_{1}\wedge\eta_{1}%
-\omega_{2}\wedge\eta_{2},$ & $d\eta_{2}=\omega_{1}\wedge\eta_{2}+\omega
_{2}\wedge\eta_{1},$ &  & \\
$\varepsilon=0,1$ & $d\eta_{3}=\varepsilon\omega_{1}\wedge\omega_{2},$ &
$d\eta_{4}=\omega_{1}\wedge\eta_{3}+\alpha\omega_{2}\wedge\eta_{3},$ &
$d\omega_{1}=d\omega_{2}=0\ $ & \\
$N_{6,24}$ & $d\eta_{1}=\omega_{1}\wedge\omega_{2},$ & $d\eta_{2}=\omega
_{2}\wedge\eta_{2},$ & $\ $ & \\
& $d\eta_{3}=\omega_{1}\wedge\eta_{3},$ & $d\eta_{4}=\omega_{1}\wedge\eta
_{3}+\omega_{1}\wedge\eta_{4},$ & $d\omega_{1}=d\omega_{2}=0$ & \\
$N_{6,25}^{\alpha\beta}$ & $d\eta_{1}=\omega_{1}\wedge\omega_{2},$ &
$d\eta_{2}=\alpha\omega_{1}\wedge\eta_{2}+\beta\omega_{2}\wedge\eta_{1},$ &
$\ $ & \\
$\alpha^{2}+\beta^{2}\neq0$ & $d\eta_{3}=-\omega_{1}\wedge\eta_{4}+\omega
_{2}\wedge\eta_{3},$ & $d\eta_{4}=\omega_{1}\wedge\eta_{3}+\omega_{2}%
\wedge\eta_{4},$ & $d\omega_{1}=d\omega_{2}=0\ $ & \\
$N_{6,26}^{\alpha}$ & $d\eta_{1}=\omega_{1}\wedge\omega_{2},$ & $d\eta
_{2}=\omega_{2}\wedge\eta_{2},$ &  & \\
& $d\eta_{3}=\alpha\omega_{1}\wedge\eta_{3}-\omega_{1}\wedge\eta_{4},$ &
$d\eta_{4}=\omega_{1}\wedge\eta_{3}+\alpha\omega_{1}\wedge\eta_{4},$ &
$d\omega_{1}=d\omega_{2}=0\ $ & \\
$N_{6,27}^{\varepsilon}$ & $d\eta_{1}=\varepsilon\omega_{1}\wedge\omega_{2},$%
& $d\eta_{2}=\omega_{1}\wedge\eta_{1},$ &  & \\
$\varepsilon=0,1$ & $d\eta_{3}=\omega_{2}\wedge\eta_{3}-\omega_{1}\wedge
\eta_{4},$ & $d\eta_{4}=\omega_{1}\wedge\eta_{3}+\omega_{2}\wedge\eta_{4},$ &
$d\omega_{1}=d\omega_{2}=0$ & \\
$N_{6,28}$ & $d\eta_{1}=\eta_{2}\wedge\eta_{4}+\omega_{1}\wedge\eta_{1},$ &
$d\eta_{2}=\eta_{3}\wedge\eta_{4}+\omega_{2}\wedge\eta_{2},$ &  & \\
& $d\eta_{3}=-\omega_{1}\wedge\eta_{3}+2\omega_{2}\wedge\eta_{3},$ &
$d\eta_{4}=\omega_{1}\wedge\eta_{4}-\omega_{2}\wedge\eta_{4},$ & $d\omega
_{1}=d\omega_{2}=0$ & \\
$N_{6,29}^{\alpha\beta}$ & $d\eta_{1}=\eta_{2}\wedge\eta_{3}+\omega_{1}%
\wedge\eta_{1}+\omega_{2}\wedge\eta_{1},$ & $d\eta_{2}=\omega_{1}\wedge
\eta_{2},$ &  & \\
$\alpha^{2}+\beta^{2}\neq0$ & $d\eta_{3}=\omega_{2}\wedge\eta_{3},$ &
$d\eta_{4}=\alpha\omega_{1}\wedge\eta_{4}+\beta\omega_{2}\wedge\eta_{4},$ &
$d\omega_{1}=d\omega_{2}=0$ & \\
$N_{6,30}^{\alpha}$ & $d\eta_{1}=\eta_{2}\wedge\eta_{3}+2\omega_{1}\wedge
\eta_{1},$ & $d\eta_{2}=\omega_{1}\wedge\eta_{2},$ &  & \\
& $d\eta_{3}=\omega_{1}\wedge\eta_{3}+\omega_{2}\wedge\eta_{2},$ & $d\eta
_{4}=\alpha\omega_{1}\wedge\eta_{4}+\omega_{2}\wedge\eta_{4},$ & $d\omega
_{1}=d\omega_{2}=0$ & \\
$N_{6,31}$ & $d\eta_{1}=\eta_{2}\wedge\eta_{3}+\omega_{2}\wedge\eta_{1}%
+\omega_{2}\wedge\eta_{4},$ & $d\eta_{2}=\omega_{1}\wedge\eta_{2},$ &  & \\
& $d\eta_{3}=\omega_{2}\wedge\eta_{3}-\omega_{1}\wedge\eta_{3},$ & $d\eta
_{4}=\omega_{2}\wedge\eta_{4},$ & $d\omega_{1}=d\omega_{2}=0$ & \\
$N_{6,32}^{\alpha}$ & $d\eta_{1}=\eta_{2}\wedge\eta_{3}+\omega_{1}\wedge
\eta_{4}+\omega_{2}\wedge\eta_{1},$ & $d\eta_{2}=\omega_{1}\wedge\eta
_{2}+\alpha\omega_{2}\wedge\eta_{2},$ &  & \\
& $d\eta_{3}=-\omega_{1}\wedge\eta_{3}+\left(  1-\alpha\right)  \omega
_{2}\wedge\eta_{3},$ & $d\eta_{4}=\omega_{2}\wedge\eta_{4},$ & $d\omega
_{1}=d\omega_{2}=0$ & \\
$N_{6,33}$ & $d\eta_{1}=\eta_{2}\wedge\eta_{3}+\omega_{1}\wedge\eta_{1}%
+\omega_{2}\wedge\eta_{1},$ & $d\eta_{1}=\omega_{1}\wedge\eta_{2},$ &  & \\
& $d\eta_{3}=\omega_{2}\wedge\eta_{3},$ & $d\eta_{4}=\omega_{2}\wedge\eta
_{3}+\omega_{2}\wedge\eta_{4},$ & $d\omega_{1}=d\omega_{2}=0$ & \\
$N_{6,34}^{\alpha}$ & $d\eta_{1}=\eta_{2}\wedge\eta_{3}+\omega_{1}\wedge
\eta_{1}+\left(  1+\alpha\right)  \omega_{2}\wedge\eta_{1},$ & $d\eta
_{2}=\omega_{1}\wedge\eta_{2}+\alpha\omega_{2}\wedge\eta_{2},$ &  & \\
& $d\eta_{3}=\omega_{2}\wedge\eta_{3},$ & $d\eta_{4}=\omega_{1}\wedge\eta
_{3}+\omega_{2}\wedge\eta_{4},$ & $d\omega_{1}=d\omega_{2}=0$ & \\
$N_{6,35}^{\alpha\beta}$ & $d\eta_{1}=\eta_{2}\wedge\eta_{3}+2\omega_{2}%
\wedge\eta_{1},$ & $d\eta_{2}=\omega_{2}\wedge\eta_{2}-\omega_{1}\wedge
\eta_{3},$ &  & \\
$\alpha^{2}+\beta^{2}\neq0$ & $d\eta_{3}=\omega_{1}\wedge\eta_{2}+\omega
_{2}\wedge\eta_{3},$ & $d\eta_{4}=\alpha\omega_{1}\wedge\eta_{4}+\beta
\omega_{2}\wedge\eta_{4},$ & $d\omega_{1}=d\omega_{2}=0$ & \\
$N_{6,36}$ & $d\eta_{1}=\eta_{2}\wedge\eta_{3}+2\omega_{2}\wedge\eta
_{1}+\omega_{2}\wedge\eta_{4},$ & $d\eta_{2}=\omega_{2}\wedge\eta_{2}%
-\omega_{1}\wedge\eta_{3},$ &  & \\
& $d\eta_{3}=\omega_{1}\wedge\eta_{2}+\omega_{2}\wedge\eta_{3},$ & $d\eta
_{4}=2\omega_{2}\wedge\eta_{4},$ & $d\omega_{1}=d\omega_{2}=0$ & \\
$N_{6,37}^{\alpha}$ & $d\eta_{1}=\eta_{2}\wedge\eta_{3}+\omega_{1}\wedge
\eta_{4}+2\omega_{2}\wedge\eta_{1},$ & $d\eta_{2}=-\omega_{1}\wedge\eta
_{3}+\omega_{2}\wedge\eta_{2}-\alpha\omega_{2}\wedge\eta_{3},$ &  & \\
& $d\eta_{3}=\omega_{1}\wedge\eta_{2}+\alpha\omega_{2}\wedge\eta_{2}%
+\omega_{2}\wedge\eta_{3},$ & $d\eta_{4}=2\omega_{2}\wedge\eta_{4},$ &
$d\omega_{1}=d\omega_{2}=0$ & \\
$N_{6,38}$ & $d\eta_{1}=\eta_{2}\wedge\eta_{3}+\omega_{1}\wedge\eta_{1}%
+\omega_{2}\wedge\eta_{1},$ & $d\eta_{2}=\omega_{1}\wedge\eta_{2},$ &  & \\
& $d\eta_{3}=\omega_{2}\wedge\eta_{3},$ & $d\eta_{4}=\omega_{1}\wedge
\omega_{2},$ & $d\omega_{1}=d\omega_{2}=0$ & \\
$N_{6,39}$ & $d\eta_{1}=\eta_{2}\wedge\eta_{3}+2\omega_{2}\wedge\eta_{1},$ &
$d\eta_{2}=\omega_{2}\wedge\eta_{2}-\omega_{1}\wedge\eta_{3},$ &  & \\
& $d\eta_{3}=\omega_{1}\wedge\eta_{2}+\omega_{2}\wedge\eta_{3},$ & $d\eta
_{4}=\omega_{1}\wedge\omega_{2},$ & $d\omega_{1}=d\omega_{2}=0$ & \\
$N_{6,40}$ & $d\eta_{1}=\eta_{2}\wedge\eta_{3}+\omega_{1}\wedge\omega_{2},$ &
$d\eta_{2}=-\omega_{1}\wedge\eta_{3},$ &  & \\
& $d\eta_{3}=\omega_{1}\wedge\eta_{2},$ & $d\eta_{4}=\omega_{2}\wedge\eta
_{4},$ & $d\omega_{1}=d\omega_{2}=0$ &\\\hline
\end{tabular}
\end{table*}

\bigskip

\end{document}